\documentclass[12pt]{article}
\usepackage[utf8]{inputenc}
\usepackage[margin=1in]{geometry}
\usepackage{hyperref}
\usepackage{amsmath}

\title{What Remains Human in Mathematics in the Age of AI} 

\vspace{1cm}

\author{ 
    Amir Moradifam\footnote{Department of Mathematics, University of California, Riverside, California, USA. E-mail: amirm@ucr.edu. }
}
\date{}

\begin{document}

\maketitle

\vspace{1cm}
For centuries, learning mathematics has meant learning to struggle productively. A student sits with a problem, tries an idea, hits a wall, tries again, and slowly discovers what matters and what does not. Even when the topic is elementary, that process is doing real work. It builds the habit of keeping track of assumptions, the ability to move between concrete examples and general principles, and the confidence to continue thinking when the answer is not immediate. Mathematics has long been valued beyond its applications precisely because of this role, as it cultivates creativity, patience, attention, judgment, and problem-solving skills.

Artificial intelligence changes this picture in a way that goes beyond the introduction of another computational tool. For the first time at scale, it becomes possible to produce mathematically correct work without necessarily engaging in the process that traditionally generates understanding. A student faced with a problem can now obtain a complete solution, often accompanied by a fluent explanation, within seconds. The resulting work may appear indistinguishable from genuine understanding, even when the underlying reasoning has never been developed.

Mathematics has repeatedly evolved alongside tools that changed the balance between routine calculation and conceptual thought, from symbolic notation and logarithm tables to calculators and computer algebra systems. The relationship between mathematics and computation has therefore always been intimate. For centuries, mathematics shaped the development of algorithms and formal reasoning, while computers gradually extended our ability to calculate and simulate. The recent emergence of large-scale artificial intelligence systems, however, marks a qualitative shift in that relationship. AI systems no longer merely execute procedures. They increasingly participate in activities associated with mathematical work itself, including formulating conjectures, manipulating symbolic expressions, generating proofs and counterexamples, and assisting in the communication of ideas.

This development necessitates a rigorous examination of the epistemology of AI-generated mathematics. Historically, the discipline has grappled with computational opacity, a concept brought to the forefront by the computer-assisted proof of the four-color theorem \cite{tymoczko1979}. In such cases, the algorithmic steps are transparent in principle, but there are simply too many for a human to verify manually. Modern AI systems, particularly large language models, introduce a different form of epistemic opacity. The internal process by which a neural network suggests a lemma or frames an argument is not directly interpretable to the user. Furthermore, general-purpose language models can still struggle with the global structure of mathematical arguments. Their autoregressive generation supports strong local coherence but does not by itself guarantee a globally consistent deductive arc across a long proof. While larger context windows mitigate this challenge to some degree, and human working memory is likewise limited, mathematicians compensate by constructing rich conceptual models, external representations, and strategic abstractions that guide long arguments. Although modern language models employ increasingly sophisticated memory and reasoning mechanisms, most do not operate on an explicit, mechanically verifiable proof state of the kind maintained by formal proof assistants. Formal proof assistants achieve their reliability through symbolic representations that are continuously checked against formal logical rules, rather than through conceptual understanding.

This raises a natural question about the practice of mathematics when machines can now engage, however imperfectly, in activities long regarded as distinctly human. Because mathematics serves as the foundational language for the sciences, its ongoing transformation offers a critical case study for the broader societal impact of AI on human knowledge systems. The impact reaches into three intertwined areas: mathematics education, the research enterprise, and the role mathematics has traditionally played in shaping how people think.

The significance of this moment lies not only in the power of the technology, but in the subtle way it alters the culture of mathematical work. Mathematics has long been defined by a particular blend of rigor, creativity, and individual insight. AI challenges us to reconsider how these qualities arise, how they are taught, and how they might evolve when part of the intellectual process is shared with computational systems. The central claim of this article is that in an era when routine tasks become easier to automate, the most human aspects of mathematics, including conceptual insight, problem framing, and creative judgment, become more visible and more valuable.

\section{Mathematics education in the age of AI}

Increasingly, a central challenge in mathematics education is no longer access to methods, but sustaining students' willingness to engage in thinking. What is at risk is not simply academic integrity, but the connection between effort and comprehension. When the act of solving is replaced by the act of requesting, there is a real danger that students come to equate exposure to a solution with ownership of an idea. Mathematics has always required more than reading. It requires attempting, failing, and reorganizing one's thinking. AI does not remove that need, but it makes it easier to bypass it. In the literature of cognitive science, this is known as cognitive offloading, which is the tendency to delegate mental effort to an external tool \cite{risko2016}. When students offload the struggle to an AI, they risk falling prey to the illusion of explanatory depth \cite{rozenblit2002}, confusing the fluent generation of a solution by the system with their own internalized understanding.

At the same time, it would be a mistake to describe AI only as a threat. Used well, it can be an unusually patient tutor. A student can ask for another example, a different explanation, a simpler starting point, or a hint rather than a full solution. Students who are reluctant to ask basic questions in a crowded classroom may find it easier to ask them privately. A motivated learner can use AI to practice, check a computation, test an approach, or explore variations of a problem. In that sense, it can lower barriers and widen access, especially for students who do not have strong support outside class.

There is also an equity angle here that is worth stating plainly. The possibility of outsourcing problem solving did not arrive with AI. Students with resources have long been able to hire private tutors, buy solution manuals, or pay for professional help that quietly erases the struggle from homework. What AI changes is not the existence of that option, but its distribution. It makes high-quality assistance cheaper and more widely available, which can reduce socioeconomic inequality in access to academic support, even as it raises new challenges for assessment and academic norms.

The same point becomes even sharper at higher levels. For a long time, access to help in advanced mathematics depended on being in the right place, having the right mentor, or finding a rare tutor who could work at the graduate level. Many topics simply do not come with an accessible support system unless one is already inside a strong department or research group. AI does not replace a good advisor, but it does make it possible to engage seriously with sophisticated material without waiting for the perfect instructor. For someone already mathematically mature, that kind of support can turn a deferred goal into real progress.

The difficulty is that the same tool can support learning or replace it, depending on how it is used. A calculator does not undermine arithmetic when students already understand place value and basic operations. It becomes harmful when it arrives before understanding. AI is similar, but more powerful, because it can provide not only answers but also plausible explanations, delivered in a confident tone even when they are wrong. That combination makes it easy for students to outsource not just computation but also reasoning. It also makes it easy for instructors to underestimate how much is being outsourced, since the written work can appear polished.

The response need not be panic or policing. It can instead be a redesign of what we ask students to do. Assignments can emphasize interpretation rather than execution, such as comparing approaches, diagnosing flawed arguments, explaining why a method works, or justifying a modeling choice. Similar concerns have motivated recent discussions on adapting mathematics education to technology-rich environments \cite{opesemowo2024}.

Yet this pedagogical shift also faces a moving target. Generative AI is increasingly capable of comparing approaches, diagnosing flaws, and producing fluent explanations. If students can simply prompt a system to generate such analyses, assessment once again risks measuring the machine's output rather than the student's understanding. The challenge is therefore not merely to design new types of assignments, but to develop forms of evaluation that require students to demonstrate ownership of the underlying reasoning. This may involve assessing how students critique, refine, and defend arguments, rather than simply evaluating the final written product.

This evolution in what we teach inevitably raises difficult questions about how we assess it. Evaluating genuine reasoning, mathematical judgment, and the ability to tolerate uncertainty is far harder to standardize than grading a final answer. Oral examinations can probe conceptual understanding through follow-up questions and guided discussion, but they are resource-intensive and difficult to implement consistently across large classes or different institutions.

One possible avenue worth exploring is the use of AI itself as an interactive assessor. Educational AI may eventually help administer oral examinations at scale. Such systems could engage students in Socratic dialogue, asking adaptive follow-up questions, probing for misconceptions, and helping instructors assess the depth of a student's reasoning without the logistical burden of human examiners. Of course, relying on AI for major assessments introduces new challenges, including ensuring the mathematical reliability of the evaluator, mitigating bias, and preventing adversarial prompting. Whether such systems ultimately prove reliable enough for high-stakes assessment remains an open question. Nevertheless, they illustrate how AI itself may become part of the solution as well as part of the challenge. Ultimately, designing scalable and reliable assessments that measure conceptual understanding rather than routine execution remains one of the most pressing practical challenges facing mathematics education.

There is also a deeper educational question underneath all of this, one with profound implications for the future scientific workforce. If mathematics has been one of the main ways we teach people to think carefully, what happens when the most visible part of mathematical work, producing solutions, can be done instantly by a machine? How we adapt mathematical education to this new paradigm of human--computer interaction will ultimately shape the cognitive habits of future scientists and the broader societal structures of knowledge production. In many fields, routine competence is still necessary, but it is no longer enough to stand out. When standard tasks can be automated or cheaply assisted, their relative value in the job market may decline. What begins to matter more is the ability to frame a problem, make good choices, notice what is missing, and push an idea beyond a template. Creativity, judgment, and the ability to use tools well become more valuable intellectual skills.

This is exactly where mathematics education can regain clarity of purpose. We may need to be more explicit with students about what the point is. The goal is not just learning techniques, but learning how to reason, tolerate uncertainty, and build an argument that can be trusted. A student who understands this is less likely to treat AI as a shortcut and more likely to treat it as an instrument. Used responsibly, AI can handle some of the mechanical burden, but the student still has to decide what to ask, how to check what comes back, how to connect it to underlying ideas, and how to judge whether an explanation is actually correct. These are not side skills. They are increasingly central intellectual skills. In the age of AI, the goals of mathematical education are not outdated. They become even more valuable, precisely because they cannot be replaced by fluent output.

\section{AI and mathematical research}
If the influence of AI in education raises questions about how mathematics is learned, its influence in research raises questions about how mathematics is created. Research has never been a purely linear process. It moves through intuition, false starts, partial computations, informal reasoning, and only later settles into the polished form of a proof. Much of this work is exploratory. One tests examples, searches for patterns, reformulates statements, or asks whether a known technique might apply in a new setting. In this phase, AI can be unexpectedly useful. It can generate examples, suggest intermediate lemmas, reorganize complicated arguments, assist with symbolic manipulation, and help researchers navigate unfamiliar areas. A mathematician moving slightly outside their specialization may use AI to sketch the landscape of a topic, recall definitions, or explore possible connections before turning to the literature. In that sense, AI can lower the cost of curiosity and broaden access to mathematical exploration.

Recent developments demonstrate that AI systems are increasingly participating in the discovery process. However, the way these systems blur the boundary between human insight and computational assistance differs fundamentally depending on the mathematical task.

First, Davies et al.~\cite{davies2021} used machine learning to detect previously unnoticed relationships in knot theory and representation theory. In knot theory, the AI analyzed high-dimensional data to reveal a surprising connection between algebraic invariants and the hyperbolic geometry of knots. The AI did not write the proof. Rather, it acted as an intuition pump, providing the empirical correlation that human mathematicians subsequently formalized into a rigorous theorem. The boundary between human and machine was blurred by outsourcing the discovery of structural correlations.

In contrast, Fawzi et al.~\cite{fawzi2022} developed AlphaTensor, which framed the search for matrix multiplication algorithms as a game played in a finite combinatorial space. By applying reinforcement learning, the AI discovered novel and faster algorithms for multiplying small matrices that outperformed some of the best human-designed algorithms. Here, the AI acted as a combinatorial search engine, automating much of the discovery process in spaces too vast and unstructured for human intuition to navigate.

Romera-Paredes et al.~\cite{romera2024} introduced FunSearch, a system that paired a large language model with an evolutionary algorithm and an automated deterministic evaluator. Instead of generating formal proofs, the system generated executable code to construct mathematical objects. It found a larger cap set in dimension eight than was previously known and improved lower bounds for the cap-set problem in extremal combinatorics. In this context, the AI acted as an automated generator of extremal structures, producing valid mathematical objects that could be checked deterministically.

Similarly, Trinh et al.~\cite{trinh2024} introduced AlphaGeometry, a neuro-symbolic theorem prover that solved 25 of 30 olympiad-level geometry problems in the IMO-AG-30 benchmark. Its architecture combines a symbolic deduction engine with a neural language model trained on a large-scale synthetic dataset to propose auxiliary geometric constructions. Within this system, the language model guides the search for useful auxiliary constructions, while the symbolic engine uses them to derive the remaining proof steps.

More recently, AI has moved from specialized research tools toward systems that can tackle difficult open problems with much less human direction. For example, in early 2026, a system combining OpenAI's GPT-5.2 Pro with Harmonic's Aristotle produced a formal Lean proof resolving Erdős Problem \#728; Sothanaphan subsequently translated the formal Lean proof into a conventional, human-readable mathematical argument \cite{sothanaphan2026}. Even more strikingly, in May 2026, an internal OpenAI reasoning model disproved the longstanding planar unit distance conjecture \cite{openai2026}. Posed by Paul Erdős in 1946, the problem asks for the maximum possible number of pairs at unit distance among $n$ points in the plane \cite{erdos1946}. The model used ideas from algebraic number theory to construct an infinite family of point sets with at least $n^{1+\varepsilon}$ unit distances for some fixed $\varepsilon>0$, contradicting the conjectured $n^{1+o(1)}$ upper bound. In this case, AI took on a much larger share of the discovery process, while human mathematicians subsequently checked, digested, and presented the result \cite{openai2026,alon2026}.

At the same time, the limitations of these systems remain substantial. AI can produce arguments that appear convincing while containing subtle gaps, incorrect references, or invalid steps. Mathematical research ultimately depends on verifiability. A proof is valuable not merely because it exists, but because it can be checked and trusted by a careful reader. Because general-purpose language models do not reliably maintain a mechanically verifiable representation of the global logical structure of an argument, their output must be treated cautiously, and responsibility for correctness remains with the human researchers who present and use the result.

Yet the central work of mathematics remains deeply human. The core of research is not the mechanical production of formal arguments, but the ability to recognize which questions matter, which structures are meaningful, and which ideas reveal something genuinely new. AI may accelerate parts of mathematical work, but, as Thurston emphasized, the discipline ultimately depends on creating and communicating understanding rather than merely producing formal arguments \cite{thurston1994}. The most productive mathematicians in the coming years may not be those who avoid AI, but those who learn to use it without allowing it to substitute for understanding.

\section{The future of mathematics}
If AI is changing how mathematics is learned and how it is practiced, the deeper question is what kind of mathematicians the future will require. For most of its history, the profession has rewarded technical mastery, persistence, and the ability to carry long chains of reasoning without error. Those qualities will not disappear, but the balance among them may shift. As routine symbolic manipulation, example generation, and even parts of exposition become easier to outsource, the value of mathematical work may move more toward interpretation, synthesis, and judgment. The mathematician of the future may spend less time executing procedures and more time deciding which problems are worth attention, which ideas are meaningful, and how different pieces of mathematics fit together.

This shift does not diminish the discipline. If anything, it makes its intellectual core more visible. Mathematics has never been only about producing proofs. It has been about seeing structure where others see complexity, asking the right question rather than merely answering one that has already been posed, and building concepts that make later discoveries possible. Mathematics has also long been valued for its remarkable ability to reveal hidden structure in the world, what Wigner famously described as its unreasonable effectiveness \cite{wigner1960}. Recent discussions have similarly suggested that the rise of AI may ultimately clarify rather than diminish the human role in mathematics, shifting attention away from routine execution and toward creativity, judgment, and deeper forms of mathematical understanding \cite{alper2026}. As routine tasks become easier to automate, asking good questions may become as important as answering them. AI may help identify promising patterns and questions, but it does not by itself determine what the mathematical community should regard as interesting, deep, or worth understanding. Those judgments remain human, shaped by taste, experience, and the collective culture of the field.

There is also a broader opportunity here. If some of the technical barriers to entry are lowered, mathematics may become more porous across its own boundaries. Researchers may move more freely between areas, revisit ideas they once abandoned, or engage with fields that previously felt too distant. Students may arrive in graduate school with a wider range of tools and curiosities. Collaboration may become more common, not only between people but between people and systems that help manage complexity. In that sense, AI may not shrink the space of mathematics, but expand the ways in which mathematicians can explore it.

At the same time, the presence of AI forces the discipline to articulate more clearly what it values. If polished arguments and fluent explanations become easier to generate, then depth, originality, and genuine understanding become more important as markers of contribution. The community may place greater emphasis on conceptual insight, on the ability to connect ideas across areas, and on the skill of recognizing when a result changes how we think rather than merely adding another statement to the literature. These qualities remain difficult to automate, and they may come to define mathematical excellence even more sharply than before.

However, if the trajectory of AI-assisted mathematics means that human mathematicians must increasingly focus on interpreting and validating machine-produced results, we are forced to confront deep epistemological questions. Will human understanding always succeed, given the sheer abundance and complexity of AI-generated theorems? Mathematics already contains partial precedents for this condition, as illustrated by the classification of finite simple groups and the four-color theorem. We accept these results and use them as tools because they have undergone extensive communal verification, even when no single mathematician is expected to hold the entire proof in mind. As AI systems produce increasingly complex combinatorial proofs, we may likewise come to accept mathematical truths whose full deductive structure lies beyond the grasp of any individual mathematician, shifting mathematics from a discipline of complete deductive transparency toward one in which verified but partially opaque results are used instrumentally.

Furthermore, what is the nature of this understanding if it occurs after a nonhuman agent provides the result? Mathematical understanding has often been treated as constructive: it deepens as a mathematician mentally reconstructs a proof. When an AI provides a result through an epistemically opaque process, human understanding may shift from being constructive toward being more interpretive or hermeneutic. The mathematician becomes less of a builder and more of an interpreter, trying to extract semantic meaning from a proof generated by an opaque computational process.

The relationship between mathematics and AI is unlikely to be one of replacement. It is more likely to become a form of partnership, uneven and sometimes uneasy, but potentially fruitful. Mathematics has always evolved with its tools, and each generation has had to decide how to use new capabilities without losing the habits of thought that define the discipline. The current moment is another such transition. If mathematicians approach AI with curiosity rather than fear, and with discipline rather than complacency, it may become not a threat to mathematical thinking but another instrument through which that thinking grows.

In the end, the question is not whether AI will change mathematics. It already has. The more interesting question is what these changes reveal about the discipline itself. Mathematics has endured for centuries not because of the tools used to practice it, but because of the habits of mind it cultivates and the clarity it demands. If AI alters the surface of mathematical work, it may also sharpen our sense of what lies at its core. What remains human in mathematics is not merely the ability to compute or even to prove, but the capacity to judge, to imagine, and to understand. Those qualities are not displaced by new tools. They are illuminated by them. \\

{\bf Acknowledgments.} The author’s work was supported in part by the National Science Foundation under grant DMS-1953620. The author also gratefully acknowledges support from a gift fund established by Turing AI.

\end{document}